\documentclass[10pt,a4paper]{article}
\usepackage[english]{babel}
\usepackage[latin1]{inputenc}
\usepackage[T1]{fontenc}
\usepackage{graphicx}
\usepackage{epsfig}
\usepackage{amsmath, amssymb}
\usepackage{color}
\usepackage{fancyhdr}
\usepackage[hypertex,dvips,colorlinks]{hyperref}
\usepackage{thumbpdf}
\usepackage{setspace}
\usepackage{vmargin}
\setmarginsrb{3cm}{0.2cm}{3cm}{2.5cm}{1cm}{1cm}{2cm}{1cm}
\definecolor{rosso}{rgb}{1,0,0}
\definecolor{nero}{cmyk}{0,0,0,1}
\definecolor{blu}{rgb}{0,0,1}

\title{Errors Theory using Dirichlet Forms, Linear Partial
  Differential Equations and Wavelets}

\author{\linespread{1} Simone Scotti\footnote{I am grateful to Nicolas Bouleau for his helpful comments, work done partially during a visit to Institut Mittag-Leffler (Djursholm, Sweden). \newline
 email: simone.scotti@unito.it 
\newline
 address: Via Real Collegio 30, 10024 Moncalieri (TO) ITALY}
\\ \\ {\it \small  Università di Torino, Ecole des
   Ponts - CERMICS} \\
  {\it \small and Collegio Carlo Alberto} }

\newtheorem{teorema}{Theorem}

\newtheorem{ipotesi}{Hypothesis}
\newtheorem{proposizione}{Proposition}
\newtheorem{remarque}{Remark}
\newtheorem{definizione}{Definition}

\newtheorem{notazione}{Notation}
\newtheorem{example}{Example}

\begin{document}

\color{nero}

\maketitle

\begin{abstract}
  We present an application of error theory using Dirichlet Forms in
  linear
  partial differential equations (LPDE). We study the transmission of an uncertainty on the
  terminal condition to the solution of the LPDE thanks to the decomposition of the
  solution on a wavelets basis.

  We analyze the basic properties and a particular class of LPDE where the wavelets
  bases show their powerful, the combination of error theory and wavelets basis justifies some
  hypotheses, helpful to simplify the computation.

\end{abstract}

\section{Introduction}

The study of the sensibility of the solution of a Partial
Differential Equation to a small perturbation of the starting
condition is a key problem in applied mathematics. The classical
approach is to define a basis of  ``perturbation functions'' and to
evaluate the variation of the PDE's solution alongs each
"direction"; this method is known as Gateaux derivative of the
solution. But this approach assumes that the perturbation of the
starting condition is deterministic, on the other hand the starting
condition is often estimated by means of some measurements, i.e.,
with mathematical language, make a statistic; therefore the result
of this estimation is a random variable, with a few known
moments\footnote{Generally we know the mean, the variance and,
maybe, the skewness and the kurtosis.}.

The probabilistic nature of the uncertainty on each estimation
pushes us to give up the Gateaux derivative approach; we follow the
methodology suggested by Bouleau, see \cite{bib:Bouleau-erreur},
this idea yields a representation of small perturbation coherent
with the truncated expansion of the perturbated solution by a small
random variable. If the estimation of starting conditions is good
the related uncertainty is very small, we may neglect order higher
than the second, that is we choose to work with Gaussian
distributions. The error theory using Dirichlet forms has two
operators, one describes the variance, called ``carré du champ'',
other one the bias, called the generator of the semigroup.

The analysis of the uncertainty transfer of the terminal function into 
the PDE's solution requires to specify a representation basis, we
choose to work with a wavelet basis and we prove some useful
properties in this case.

Summarizing, we propose a new approach to study the impact of
uncertainty on the solution of a Linear Partial Differential
Equation due to an random imprecision on the starting condition;
this method permits, first of all, the study of variance of the LPDE
solution, i.e. his sensitivity; secondly we can estimate a
covariance between the LPDE solution at two different points of
time-space domain.

The paper is organized as follows: In section 2, we present the
error theory using Dirichlet Forms and we model, mathematically, the
uncertainty on the starting condition. Section 3 is a survey of
wavelets theory. In section 4 we study the solution of a LPDE using
wavelet's decomposition and we present some particular cases where
the wavelet's properties play a crucial role. In section 5 we
describe the evolution equation for the operators of error theory in
the LPDE case. Section 6 is dedicated to an introduction to possible
applications in finance. The results are summarized in section 7.

\section{Error theory using Dirichlet forms}

To make the paper self-contained, we recall the principal ideas of
error theory using Dirichlet forms, for reference see Bouleau
\cite{bib:Bouleau-erreur}.

The first studies of error propagation, owing to Legendre, Laplace
and Gauss, date back at the beginning of the 19th century, show the
probabilistic nature of errors. But the classical approach in the
study of sensibility is to evaluate the derivative in the sense of
Gateaux with respect to the erroneous parameters. Clearly this
approach is coherent if the nature of incertitude is determinist
(i.e. the ``good'' value exist) and it is incoherent with a
probabilistic nature.

The error theory using Dirichlet forms permits to study the
sensibility with respect to a parameter, characterized by a random
uncertainty.

In order to introduce this theory, we analyze  a particular case, we
start with a parameter $\sigma_0$ and we want evaluate
$F(\sigma_0)$, where $F$ is a nonlinear regular function. If the
quantity $\sigma_0$ has a small centered error, we can model this
parameter with a random variable $\sigma = \sigma_0 + \epsilon
A[\sigma] + \sqrt{\epsilon \gamma[\sigma]} Y$, where $\epsilon$ is a
small constant that defines the order of magnitude of the
perturbation, $A[\sigma]$ represents the renormalized bias of the
random variable with respect to the parameter value $\sigma_0$,
$\gamma[\sigma]$ the variance and $Y$ is a Gaussian standard random
variable.

If we study the bias and the variance of $F(\sigma)$, use of
Taylor's formula show

\begin{equation*}
  \left\{
    \begin{array}{lll}
      \mathbb{E}[F(\sigma)-F(\sigma_0)] & = & \epsilon \, F'(\sigma_0)
\, A[\sigma] + \frac{1}{2} \, \epsilon \, F^{\prime \prime}(\sigma_0) 
\,\gamma[\sigma] + o(\epsilon) \\
      \mathbb{E}[\left\{F(\sigma)-F(\sigma_0) \right\}^2] & = & \epsilon  
\, \left[F^{\prime}(\sigma_0) \right]^2 \, \gamma[\sigma] + o(\epsilon)
    \end{array}
  \right.
\end{equation*}

We emphasize the following remark.
\begin{remarque}
  The crucial difference between the deterministic approach and the probability one
  is in the bias term that depends on the variance of the random variable $\sigma$ if
  the function $F$ is nonlinear; therefore this effect is said
  stochastic.
\end{remarque}

The original idea of error theory using Dirichlet Forms is to stop
the Taylor expansion at the first order in $\epsilon$ and to forget
the higher orders; and after search for a theory with two operators
with the following chain rule:

\begin{equation}\label{chain-rule}
  \left\{
    \begin{array}{lll}
     \mathcal{A}[F(\sigma)] & = & F'(\sigma) \, \mathcal{A}[\sigma] +
      \frac{1}{2} \, F^{\prime \prime}(\sigma) \, \Gamma[\sigma] \\
      \Gamma[F(\sigma)] & = & \left[F^{\prime}(\sigma) \right]^2 \, 
\Gamma[\sigma]
    \end{array}
  \right.
\end{equation}
where $\Gamma$ plays the role of a variance-covariance operator and
$\mathcal{A}$  a bias operator.

This theory exists and is the theory of semi-group, the operator
$\mathcal{A}$ is the generator of the semi-group and the operator
$\Gamma$ is the carré du champ associated with the Dirichlet Form of
the semi-group, see Albeverio \cite{bib:Albeverio} for a general
introduction.

We begin the introduction of error theory using Dirichlet Forms with
the main definition.

\begin{definizione}[Error Structure]
An error structure is a term $\displaystyle \left(\Omega, \, \mathcal{F}, \, \mathbb{P}, \,
  \mathbb{D}, \, \Gamma \right)$ where:

\begin{enumerate}
\item $\displaystyle \left(\Omega, \, \mathcal{F}, \, \mathbb{P}\right)$  is a
  probability space;
\item $\mathbb{D}$ is a dense sub-vector space of $\displaystyle L^2 \left(\Omega, \,
    \mathcal{F}, \, \mathbb{P}\right)$;
\item $\Gamma$ is a positive symmetric bilinear application  form $\mathbb{D} \times
  \mathbb{D}$ into $L^1 \left(\Omega, \,
    \mathcal{F}, \, \mathbb{P}\right)$, such that if $F$ and $G$ are two functions
  belonging $C^1$ and Lipschitzian class, $u$ and $v \in \mathbb{D}$, then we have
  $F(u)$ and $G(v) \in \mathbb{D}$ and

  $$
  \Gamma[F(u), \, G(v)] = F^{\prime}(u) G^{\prime}(v) \Gamma[u, \, v] \; \;
  \mathbb{P} \, a.s;
  $$
\item The space $\mathbb{D}$ equipped with the norm $ \displaystyle
  ||u||_{\mathbb{D}}= \left(||u||_{L^2} +  \mathbb{E}\left[\Gamma[u, \,u] \right]
  \right)$ is complete.

\end{enumerate}

\end{definizione}

When the two arguments of gamma operator are the same, we can
simplify the notation $\Gamma [u] \equiv \Gamma [u, \, u]$.

The knowledge of the carré du champ $\Gamma$ and the probability
$\mathbb{P}$ permits to associate a unique strongly continuous
contraction semi-group $\left(P_t \right)_{t \geq 0}$ thanks to
Hille Yosida theorem (see Albeverio pages 20-26). This semi-group
has a generator $(\mathcal{A})$ defined on a domain
$\mathcal{D}\mathcal{A}$ subspace  of $\mathbb{D}$ and it is a
self-adjoint operator that satisfies the chain rule (\ref{chain-rule})
when $F \in C^2$, $\sigma \in  \mathcal{D}\mathcal{A}$ and
$\gamma[\sigma] \in L^2 (\mathbb{P})$

The operator $\Gamma$ is bilinear, like the variance-covariance one,
that makes computations awkward to perform. It is possible to
surmount this difficulty  by means of a new operator, the sharp,
which, in some sense, is a square root version of $\Gamma$, see
Bouleau \cite{bib:Bouleau-erreur} page78-80.

\begin{definizione}[Sharp operator]\label{def:sharpe}
Let $\displaystyle \left(\Omega, \, \mathcal{F}, \, \mathbb{P}, \,
  \mathbb{D}, \, \Gamma \right)$ be and error structure and $\displaystyle
\left(\widehat{\Omega}, \, \widehat{\mathcal{F}}, \, \widehat{\mathbb{P}} \right)$ be
a copy of the probability space. If the space $\mathbb{D}$ is separable, there exist
an operator sharp, denoted $()^{\#}$, with these three properties:

\begin{enumerate}
\item $\forall u \in \mathbb{D}$, $u^{\#} \in L^2(\mathbb{P} \times
  \widehat{\mathbb{P}})$;
\item $\forall u \in \mathbb{D}$, $\Gamma[u, \, u] =
  \mathbb{E}^{\widehat{\mathbb{P}}}\left[\left(u^{\#} \right)^2 \right]$;
  \item $\forall u \in \mathbb{D}$ and $F \in C^1$, $\left(F(u)\right)^{\#} =
    F^{\prime}(u) \, u^{\#}$.
\end{enumerate}

\end{definizione}

We illustrate the definition of error structure and related topics by means of an
example.

\begin{example}[real space]

  The simplest structure on the real space is the Orstein-Uhlenbeck
  structure:

  $$
  \left(\Omega, \, \mathcal{F}, \, \mathbb{P}, \,
    \mathbb{D}, \, \Gamma \right) =  \left(\mathbb{R}, \, \mathcal{B}(\mathbb{R}), \,
    \mu, H^1(\mu), \Gamma[u, \,u] = \left\{u^{\prime}\right\}^2 \right)
  $$
where $\mathcal{B}(\mathbb{R})$ is the Borel $\sigma$-field of $\mathbb{R}$, $\mu$ is a gaussian measure and $H^1(\mu)$ is the first Sobolev space with respect
to the measure $\mu$, i.e. $u \in H^1(\mu)$ if $u \in L^2(\mu)$ and $u^{\prime}$ in
the distribution sense belongs to $L^2(\mu)$.

The associate generator has the following domain:

$$
\mathcal{D}\mathcal{A} = \left\{u \in L^2(\mu): \, u^{\prime \prime} - x \,f^{\prime}
\text{ in the distribution sense belongs to } L^2(\mu)\right\}
$$

and the generator operator is

$$
\mathcal{A}[u] = \frac{1}{2} u^{\prime \prime} - \frac{1}{2} I \cdot u^{\prime}
$$

where $I$ is the identity map on $\mathbb{R}$

The associate sharp operator is

$$
(u)^{\#}(x, \, y) =  u^{\prime}(x) G(y)
$$

where $G$ is a fixed function such that

$$
\int_{\mathbb{R}}  G^2(y) \, d\mu(y) = 1
$$

\end{example}

\subsection{Error structure on functional spaces}\label{sec:error-function}

The idea of error structure also can be used when the probability
space $\Omega$ has infinite dimension, in this case some constraints
must be verified see Bouleau \cite{bib:Bouleau-erreur} pages 59-64.
In particular, afterwards in this paper, we need to define an error
structure on the space of functions $L^2$.

Let $\Omega$ be a space of functions, we assume that it exists a
basis (perhaps numerable) $\{\psi_i\}_{i \in \mathbb{N}}$ of this
space; therefore if $f \in \Omega $ then

$$
f(x) = \sum_{i=0}^\infty a_i \psi_i(x)
$$

The simplest way in order to define an error structure on the space
$\Omega$ is to randomize the coefficients $a_i$, to define an error
structure on each subspace generated by each element of the basis
and take the product structure as the error structure on the
functional space $\Omega$. Bouleau has proved that an infinite
product of error structures is still one, under some constraints on
the $\Gamma$ operator, see \cite{bib:Bouleau-erreur} page 62.

An easy way to force the wellness of the product of error structures
is to truncate the previous sum, see Bouleau
\cite{bib:Bouleau-erreur} page 84 for a practical example.

On each factor we have the choice on the structure of gamma
operator; a first idea is to suppose that the uncertainties on two
different factors are uncorrelated:

\begin{ipotesi}[independence]\label{ipo:independance}

  The operator $\Gamma$ verifies:

\begin{equation}
\Gamma[a_m, \, a_n] = 0 \text{ }  \text{ when } \text{ } m \neq n
\end{equation}

\end{ipotesi}

This hypothesis simplifies the course of the analysis of uncertainty
diffusion, we emphasize the nature of this hypothesis in the
following remark.

\begin{remarque}
  It is clear, first, that the hypothesis of uncorrelation must be have an external
  justification, the error theory using Dirichlet forms cannot force a particular
  structure on the error correlation. Secondly, the choice of the basis is crucial,
  since, if we start from an error structure, associated at a basis, that verifies
  the hypothesis (\ref{ipo:independance}), and we change the basis, the image error
  structure through the new basis  presents some correlations except in exception
  circumstances.

  Therefore, it is clever to choice a basis with some features in order to use
  these
  properties to justify the hypothesis of uncorrelation.
\end{remarque}

Lastly we have to specify the choice about the operator $\Gamma$ on
each subspace generated by each element of the basis; many options
are possible but we prefer the following

\begin{ipotesi}[proportionality]\label{ipo:proportionality}

 The operator $\Gamma$ verifies:

\begin{equation}
\Gamma[a_n, \, a_n] = \gamma(n)\, a_n^2
\end{equation}

\end{ipotesi}

The advantage of this hypothesis is to yield a constant proportional
error, i.e. each factor on the basis decomposition contributes
towards the uncertainty proportionally to their contribution at the
decomposition of the function; therefore, if a function as not
projection on a subspace, the same subspace do not contribute to the
uncertainty. The function $\gamma(n)$ permits to leave some freedom
at the class of error structure.

Considering the fact that we have an independent error structure on
each subspace generated by the basis decomposition, thus the
Mokobodzki hypothesis is verified, the error structure admits a
sharp operator, see Bouleau \cite{bib:Bouleau-erreur} pages 78-82,
and it verifies the properties:

\begin{equation}\label{eqn:sharp-indep}
\begin{array}{rcl}
  a_n^{\#} & = & \sqrt{\gamma(n)} a_n \widehat{a}_n  \\
  \widehat{\mathbb{E}}\left[\widehat{a}_n^2  \right]  & = & 1
  \end{array}
\end{equation}

\section{Wavelets theory}

In this section we give a short presentation of the theory of
wavelets, in order to preserve the self-contained of this paper,
this introduction  follows Hardle et al. \cite{bib:Hardle} and
Mallat \cite{bib:Mallat}. We fix a function, called father wavelet,
$\phi \in L^2(\mathbb{R})$, such that the family $\{\phi_{0,\, k} =
\phi(\cdot - k), k\in \mathbb{Z}\}$ is an orthonormal
system\footnote{The easy way is to consider a compactly supported
$\phi$.}.

We define the linear space (sub-space of $L^2$)
\begin{displaymath}
  V_0 = \left\{f(\cdot) \,\left| \;  f(\cdot) = \sum_{k=-\infty}^{\infty} c_k
      \phi_{0, \, k}(\cdot) \right. \right\}.
\end{displaymath}

From this original space we can define a chain of sub-spaces $\{V_i \}_{i \in
  \mathbb{Z}}$ by the relation

\begin{displaymath}
  f(\cdot) \in V_i \text{ } \text{  iff } \text{ } f(2 \cdot) \in V_{i-1}
\end{displaymath}

these spaces are called ``generated'' by the function $\phi$. Mallat
and Meyer introduce the concept of multiresolution analysis in the
years 1988-1990

\begin{definizione}[MRA]
  A chain of sub-spaces $\{V_i \}_{i \in
    \mathbb{Z}}$, ``generated'' by a function $\phi$ is called a Multiresolution
  Analysis [MRA] if  the two following conditions are held:

  \begin{itemize}
  \item $V_i \subset V_{i+1} $ for all i and
  \item $\bigcup_{i} V_i$ is dense in $L^2$
  \end{itemize}
\end{definizione}

In this case the function $\phi$ is called the ``father Wavelet''.

To define an orthogonal basis of the $L^2$ space we must define a
sequence of orthogonal spaces, since the chain $\{V_i\}$ is
decreasing sequence. Define $W_i$ the orthogonal complement of $V_i$
into $V_{i-1}$, for all i; we find  that $V_i = V_0 \oplus
\bigoplus_{j=1}^i W_j$ and, thanks to the second property of MRA, we
have:

\begin{displaymath}
  V_0 \oplus
\bigoplus_{j=1}^{\infty} W_j \text{ is dense in } L^2
   \end{displaymath}

   We can fix an orthonormal basis $\{\psi_{i \, k}\}_{k \in \mathbb{Z}}$ in each space $W_i$,
   Mallat and Meyer show that we can fix $\psi_{i, \, k}(\cdot) = \sqrt{2^{-j}} \psi( 2^{-j}
   \cdot -k)$ where $\psi(\cdot)$ is a function depending of $\phi(\cdot)$, see Mallat
   \cite{bib:Mallat}, page 233 for the explicit relation; the function $\psi(\cdot)$
   is called the ``mother Wavelet''.

   We conclude that the original function $f(x)$ has a unique representation in term
   of the following series:

   \begin{displaymath}
f(x) = \sum_k \alpha_k \phi_{0, \, k}(x) + \sum_i \sum_k \beta_{i, \, k} \psi_{i, \, k}(x)
   \end{displaymath}

\subsection{Daubechies wavelets}

The construction of the wavelet basis depends on the choice of the father wavelet;
this choice is not constrained, many options are possible, in
this section we present a particular class of wavelets called Daubechies Wavelets,
see Daubechies \cite{bib:Daubechies-art}.

The first historical wavelet is the Haar's basis, its father wavelet
is $\phi (x) = \mathbf{1}_{[0, \, 1]}$; this basis has some
advantage, in particular the boundary support and the quick
computation; unluckily the Haar mother wavelet is discontinuous,
therefore this basis misapproximates all continuous function, this
roughly approximation induces that the coefficients $\beta_{i, \,
k}$ do not decrease fast with the rescaling index $i$.

Another possibility is the Riesz's class, but the choice of the
Riesz's bases approach cause uncompactly supported father and mother
wavelets, see \cite{bib:Hardle} chapter 6, and it is clear that a
function with compact support is easy to treated numerically.

Actually the simpler known class of wavelets with continuous and
compactly supported father and mother wavelets is the Daubechies'
wavelets, see \cite{bib:Hardle} chapter 7 or \cite{bib:Mallat} pages
246-251.

The definition of Daubechies' father wavelet starts from its Laplace
transform.

\begin{displaymath}
\widehat{\phi} (\omega) = \prod_{j=0}^{\infty} m_0(2^{-j} \omega)
\end{displaymath}
where $m_0$ is a $2\pi$-periodic function, that must verify some
constraints, see \cite{bib:Hardle} chapter 5 and 6 or
\cite{bib:Mallat} pages 222-231; in particular an usual choice for
$m_0$ is a trigonometric polynomial:

  \begin{displaymath}
     m_0(\omega) = \sum_k h_k e^{-i \,k \, \omega}
    \end{displaymath}

The Daubechies' basis of order p must verify that the function
$m_0(\omega)$ has a zero of order p for $\omega = \pi$; we conclude
with the equation of the Fourier transform of mother wavelet:

\begin{displaymath}
  \widehat{\psi}(\omega) = \frac{1}{\sqrt{2}}\, \overline{m_0\left(2^{-1} \omega + \pi
    \right)} \, e^{-i \, \frac{\omega}{2}} \, \widehat{\phi}\left( \frac{\omega}{2} \right)
  \end{displaymath}

the principal advantages of these bases are:

\begin{itemize}
\item{\it regularity:} the Daubechies' father and mother wavelets are uniformally $\gamma$-Lipchitz
  where the parameter $\gamma$ grow with the order p, in particular if p is grater
  than three the wavelets are differentiable, see Daubechies et al. \cite{bib:Daubech-Laga}.

  \item{\it compact support:} the support of Daubechies' father and mother wavelets are
    length $2p -1$.
  \end{itemize}

In conclusion we have a representation of the function $f(x)$ via a
Daubechies' decomposition on wavelets:

 \begin{equation}\label{eqn:wavelets-decomposition}
f(x) = \sum_k \alpha_k \phi^D_{0, \, k}(x) + \sum_i \sum_k \beta_{i, \, k} \psi^D_{i, \, k}(x)
   \end{equation}

\section{Linear partial differential equation}

In this section we want study the interaction between wavelets and
linear partial differential equations with a terminal condition
$f(x)$. In the previous paragraph we have decomposed the function
f(x) into a Daubechies' wavelets basis; now we study the evolution
of this decomposition through a linear partial differential equation
and we focus on the properties of this evolution.

We consider a Linear Partial Differential Equation

\begin{equation}\label{eqn:LPDE-generale}
\left\{
\begin{array}{rcl}
\displaystyle \frac{\partial Q}{\partial t}(t, \, x) +
\frac{\sigma^2(x, \, t)}{2} \frac{\partial^2 Q}{\partial x^2} (t, \,
x) + \mu(x, \, t) \frac{\partial
  Q}{\partial x} & = & \displaystyle 0 \\
\scriptstyle  &  \scriptstyle  & \scriptstyle  \\
 \displaystyle Q(T, \, x)   & = & \displaystyle f(x)
\end{array}
\right.
\end{equation}

The Feynman-Kac formula tells us that the solution of LPDE
(\ref{eqn:LPDE-generale}) is

\begin{equation}\label{eqn:LPDE-solution}
Q(t, \, x ) = \mathbb{E}\left[f(X_T) | X_t = x \right]
\end{equation}

where the Ito process $X_t$ verifies the following SDE

\begin{equation}\label{eqn:LPDE-EDS}
dX_t = \mu(t, \, X_t) \, dt + \sigma(t, \, X_t) \, dW_t
\end{equation}
With $W_t$ is a Brownian motion.

Since the linearity of the PDE  (\ref{eqn:LPDE-generale}), the
terminal condition and the solution admits a decomposition on basis,
in particular a wavelets basis. Then, if the terminal condition is
written as following

\begin{equation*}
f(x) = \sum_k \alpha_k \phi^D_{0, \, k}(x) + \sum_i \sum_k \beta_{i, \, k} \psi^D_{i, \, k}(x)
\end{equation*}
we can study the evolution of each factor and prove that the solution has a similar decomposition.

\begin{notazione}
In order to simplify the notation we shift the time coordinate in such a way that we bring the final time T at zero.
Clearly each time t smaller than T became negative and the terminal condition became $Q(0,\, x)=f(x)$. 
\end{notazione}

\subsection{Diffusionlets}

In this subsection we analyze a particular case in which the
wavelets properties have a key role; a basilar study has been done
by Shen and Strang, see \cite{bib:Shen-Strang}, when the PDE is the
heat equation $U_t(t, \, x) = c \triangle U(t, \, x) $.

Shen and Strang introduce the notation of mother and father heatlets (given a choice
of wavelets basis), these are the
heat evolution of the mother and father wavelets.

In a similar way, if we consider a LPDE and a wavelets basis, we can
define an object called Diffusionlet:

\begin{definizione}[Diffusionlets]
The solution of a LPDE with a wavelet as terminal condition is
called diffusionlet.
\end{definizione}

The diffusionlets associated at a LPDE are denoted:

\begin{equation}\label{eqn:LPDE-Wavelets-solution}
\begin{array}{lll}
  \Phi^D_{0, \, k}(t, \, x) & = & \mathbb{E}\left[\phi^D_{0, \, k}(X_0) | X_t = x
  \right] \\
   \Psi^D_{i, \, k}(t, \, x) & = & \mathbb{E}\left[\psi^D_{i, \, k}(X_0) | X_t = x
  \right]
  \end{array}
\end{equation}

we can make some remarks:

\begin{remarque}
  In Daubechies' case, the father and mother wavelets are compactly supported and bounded,
  therefore the variance of the solutions $\Phi^D$ and $\Psi^D$ are, generally,
  smaller than the variance of the solution $Q(t, \, x )$
\end{remarque}

\begin{remarque}
  The Daubechies' mother wavelet $\psi^D(x)$ and its rescaling functions $\psi^D_{i,
    \, k}(x)$ have p vanishing moments (where p is the order of Daubechies's wavelet), so the
  associate diffusions depend mainly on the assymmetries and we can suppose that the
  contributions of high order wavelets vanish very quickly with the time.
\end{remarque}

\begin{remarque}
  The study of the solution of the LPDE  (\ref{eqn:LPDE-Wavelets-solution}) requires to
  solve the same LPDE with each wavelet basis as terminal value.
\end{remarque}

The least remark underline the main problem with this approach, we
have earned a high precision on the estimation of the solution but
the price to pay is that we need to solve many times the same
problem in order to have a good estimation. But we introduce a class
of LPDE where this difficulty has an easy answer.

We consider a special case of the previous LPDE, we fix the function
$\mu(t, \,x)= r x$  and the function $\sigma(t, \, x)$ to be a power
of the variable x.

The LPDE becomes

\begin{equation}\label{eqn:LPDE-diffusionlets}
\left\{
\begin{array}{rcl}
\displaystyle \frac{\partial Q}{\partial t}(t, \, x) + \frac{\sigma^2 \, x^{2\lambda}}{2}
\; \frac{\partial^2 Q}{\partial x^2} (t, \, x) + r \, x \, \frac{\partial
  Q}{\partial x} & = & 0 \\
\scriptstyle & \scriptstyle & \scriptstyle \\
\displaystyle  Q(0, \, x)   & = & \displaystyle f(x)
\end{array}
\right.
\end{equation}
and the associate diffusion is

\begin{equation}\label{eqn:LPDE-diffusionlets-EDS}
dX_t = r \, X_t \, dt + \sigma \, X_t^{\lambda} \, dW_t
\end{equation}

A change of numeraire proves that the good solution of LPDE
(\ref{eqn:LPDE-diffusionlets}) can be write as $Q(t,\, x) = \widetilde{Q}(t, \, x)
e^{-r\,t}$ where $\widetilde{Q}(t, \, x)$ verifies the following LPDE

\begin{equation}\label{eqn:LPDE-diffusionlets-discounted}
\left\{
\begin{array}{rcl}
\displaystyle \frac{\partial \widetilde{Q}}{\partial t}(t, \, x) + \frac{\sigma^2\, x^{2\lambda}}{2}
\;\frac{\partial^2  \widetilde{Q}}{\partial x^2} (t, \, x) & = & 0 \\
\scriptstyle & \scriptstyle & \scriptstyle \\
\displaystyle  \widetilde{Q}(0, \, x)   & = & \displaystyle f(x)
\end{array}
\right.
\end{equation}

\begin{remarque}[Self-similarity]\label{rem:self-similarity}
  This equation presents a rescaling invariance: if $ \widetilde{Q}(t, \, x)$ is a
solution of LPDE (\ref{eqn:LPDE-diffusionlets-discounted}) with terminal value $g(x)$,
then also $ \displaystyle
\widetilde{Q}(\alpha^{2-2 \lambda} \,t, \,
\alpha \, x)$ is a solution of LPDE with terminal value $g(\alpha \, x)$.
\end{remarque}

We can prove some properties and theorems of the wavelets diffusion
through a LPDE of type (\ref{eqn:LPDE-diffusionlets}), each results is
a generalization of an equivalent theorem on heatlets decomposition,
see Shen and Strang \cite{bib:Shen-Strang}, the crucial hypothesis
used in all theorems is the linearity of PDE.

\begin{proposizione}[Refinement]
  Let $\phi(x)$ and $\psi(x)$ be, respectively, the father and mother wavelets of a
  wavelets basis.
  Suppose that the LPDE
  (\ref{eqn:LPDE-diffusionlets-discounted}) has a unique solution for each terminal
  function $f(x)$;
  let $\Phi^D(t, \, x)$ and  $\Psi^D(t,\, x)$ be the
  solutions of the LPDE (\ref{eqn:LPDE-diffusionlets-discounted}) with terminal
  conditions, respectively, the father and the mother wavelets. If

  \begin{equation}\label{prop_1:hypo}
    \begin{array}{lll}
\phi(x) & = & \displaystyle 2 \sum_{n \in \mathbb{Z}} \, h_n \, \phi(2x -n)  \\
 \psi(x) & = & \displaystyle 2 \sum_{n \in \mathbb{Z}} \, g_n \, \phi(2x -n)
\end{array}
\end{equation}
with
\begin{displaymath}
\sum_{n \in \mathbb{Z}} (h_n^2 + g_n^2) < \infty
\end{displaymath}

Then
\begin{equation}\label{prop_1:result}
  \begin{array}{lll}
\displaystyle \Phi^D(t, \,x) & = &\displaystyle 2 \, \sum_{n \in \mathbb{Z}} h_n \, 
\Phi^D\left((2^{(2-2\lambda)} \,t, \, 2x -n \right)  \\
\scriptscriptstyle & \scriptscriptstyle &\scriptscriptstyle \\
\displaystyle \Psi^D(t, \,x) & = & \displaystyle 2 \, \sum_{n \in \mathbb{Z}} g_n \,
\Phi^D\left( 2^{(2-2\lambda)} \, t, \, 2x -n \right)
\end{array}
\end{equation}

\end{proposizione}

{\it Proof}: It is easy to check that $ \displaystyle
\Phi^D(2^{(2-2\lambda)} t, \, 2x -n)$ is the solution of LPDE
(\ref{eqn:LPDE-diffusionlets-discounted}) with terminal value
$\phi(2x-n)$, see remark \ref{rem:self-similarity}. Then the
evolution of the functions $\phi(x)$ and  $\psi(x)$ given the
left-hand sides of equations (\ref{prop_1:result}) and the
decomposition of the same functions, thanks to relations
(\ref{prop_1:hypo}) given the right-hand sides; by uniqueness the two
solutions must be equal.
\begin{flushright}
$\Box$
\end{flushright}

\begin{teorema}[Diffusionlets
decomposition]\label{teo:diffusionlets}
  Suppose that f(x) belongs to $L^2$ and then $f(x)$ admits a decomposition of type
  (\ref{eqn:wavelets-decomposition}) and suppose that LPDE
  (\ref{eqn:LPDE-diffusionlets-discounted}) admits a unique solution
  when the terminal condition is $\widetilde{Q}(0, \,x)= f(x)$
   then the solution of LPDE
  (\ref{eqn:LPDE-diffusionlets-discounted}) with terminal value $f(x)$ is given by

  \begin{equation}\label{eqn:diffusionlets-decomposition}
    \widetilde{Q}(t, \,x) = \sum_k \alpha_k  \; \Phi^D(t, \, x-k) +
    \sum_i \sum_k \beta_{i, \, k} \; \Psi^D\left(2^{(2-2\lambda)i} \,t, \, 2^j \, x  -k\right)
  \end{equation}

\end{teorema}

{\it Proof:} Since $f(x) \in L^2(\mathbb{R})$ and $\left\{\phi^D_{k}
\, \, \psi^D_{i,
    \, k}(x)\right\}$ is an orthonormal basis of $L^2(\mathbb{R})$,
    the wavelets expansion of $f(x)$ converges to $f(x)$ in
    $L^2$-norm.

It is easy to check that $\displaystyle \Psi^D\left(2^{(2-2\lambda)i}
\, t, \, 2^j \,x -k \right)$ is the evolution of $\psi^D_{i, \,
k}(x)$ and $\Phi^D(t, \, x-k)$ is the evolution of $\phi^D_{0,\,
k}(x)$, see remark \ref{rem:self-similarity}.

By uniqueness of the solution the diffusion of the wavelets
expansion of terminal condition $f(x)$ converges to the diffusion of
$f(x)$ in $L^2$-norm
\begin{flushright}
$\Box$
\end{flushright}


\begin{remarque}[Diffusionlets Advantages]\label{rem:advantage}
  The key advantage of the diffusionlets is the independence of the initial state.
  Therefore we can store the solution of the LPDE (\ref{eqn:LPDE-diffusionlets-discounted})
  with the father wavelet as terminal condition, this solution can be estimate with an high
  degree of precision due to the compactly supported and bounded father
  wavelet in Daubechies case.

  The solution of the LPDE (\ref{eqn:LPDE-diffusionlets-discounted})
  with terminal value $f(x)$
  can be evaluate using this strategy:

\begin{enumerate}
    \item study the decomposition coefficients of the function $f(x)$ into the
    wavelets basis, using FWT (fast wavelet transform) as an
    example, see Mallat \cite{bib:FWT};

    \item reconstruct the solution with $f(x)$ as terminal condition
    using the store solution for the father wavelet, the coefficients
    estimated and the result of theorem \ref{teo:diffusionlets}.
\end{enumerate}

The second, crucial,  advantage is the fact that Daubechies wavelets
of order p have p vanishing moments, this fact, combined with the
smoothing property of parabolic partial differential equation
(property true for LPDE (\ref{eqn:LPDE-diffusionlets-discounted}) far
from $x=0$) force a fast convergence of the mother wavelet to 0 with
time t; this effect is magnified by the scaling effect when
$\lambda<1$, see equation (\ref{eqn:diffusionlets-decomposition}).

\end{remarque}

\begin{remarque}[Diffusionlets disadvantage]
    Due to the linearity of the LPDE
    (\ref{eqn:LPDE-diffusionlets-discounted}), diffusionlets are not
    compactly supported when $t \neq 0$.
    This fact force the necessity to estimate an "essential
    support" for the diffusion of father and mother wavelets, i.e. the
    region where the norm of father and mother diffusionlets are bigger than
    a reference level $\epsilon$, related to the asked sensibility
    of  the searched solution.
\end{remarque}

    This "essential support" given a length, depending on time t,
    useful to determine the number of decomposition elements for each
    scaling level, i.e. the number of $k$ that it is necessary
    consider for a well-estimation of the solution before changing the coefficient
    $i$.

The smoothing effect, emphasized between the advantages, permits to
define an order $I(\epsilon)$ beyond which the remainder is smaller
than $\epsilon$

    The approximated solution becomes:

    \begin{equation}
        \widetilde{Q}(t, \, x) \simeq \sum_{k=-K_0(\epsilon)}^{K_0(\epsilon)}
        \alpha_k \; \Phi^D(t, \, x-k) + \sum_{i=0}^{I(\epsilon)}
        \sum_{k= -K_i(\epsilon)}^{K_i(\epsilon)} \beta_{i, \, k} \;
        \Psi^D\left(2^{(2-2\lambda)i}\, t, \, 2^i \, x -k\right)
    \end{equation}

We conclude with a negative remarque on the diffusionlets basis:

\begin{remarque}[non-orthogonality of diffusionlets basis]
The class of functions generated by the diffusion of a wavelets
basis is a basis of the $L^2$-subspace characterized by the
diffusion itself\footnote{The proof of this fact is very simple, we
must check only the independence of each vector of the basis, that
follows on the uniqueness of the solution.}, i.e., for each time
$s$, the smaller subspace of $L^2$ that contains the functions that
can be a solution of the LPDE at time $s$ with a terminal value
belongs to $L^2$.

However, the basis $\displaystyle \left\{\Phi^D_k(t, \, x), \;
\Psi^D_{i,\, k}(t, \, x) \right\}_{i,\, k}$, where $\displaystyle
\Phi^D_k(t, \, x) = \Phi^D(t, \, x-k)$ and $\displaystyle
\Psi^D_{i,\, k}(t, \, x) = \Psi^D\left(2^{(2-2\lambda)i}\, t, \, 2^i
\,x -k\right) $, is not, generally, orthogonal, especially owing to
the scaling factor in time.

\end{remarque}

\section{Sensibility of LPDE solution}

In this section we suppose that the terminal condition $f(x)$ of a
LPDE are erroneous and we study the diffusion of this uncertainty.
The starting point is to define an error structure on a functional
space, see subsection \ref{sec:error-function}. We use the
decomposition of $f(x)$ into a wavelets basis, see equation
(\ref{eqn:wavelets-decomposition}); and we set the coefficients
$\alpha_k$ and $\beta_{i, \, k}$ to be random.

We define an error structure on each subspace, generated by each
element of the wavelets basis, in accord with hypotheses
\ref{ipo:independance} and \ref{ipo:proportionality}, we can assume
that the error structures on each subspace are independent and the
uncertainty is proportional to the estimate parameter, i.e.
$\alpha_k$ or $\beta_{i, \, k}$ depending on the cases. Now we can
study the variance caused by the uncertainty on the terminal value.

\subsection{uncertainty on the solution}

We can prove three results:

\begin{proposizione}[Variance of terminal condition]

The terminal condition $f(x)$ has the following variance:

\begin{equation}\label{eqn:gamma-terminal-function}
  \Gamma\left[ f(x) \right]  =   \sum_k \gamma(k) \, \alpha_k^2  \,  \left[\phi^D_{0, \,
    k}(x)\right]^2 + \sum_i \sum_k \gamma(i, \, k) \, \beta^2_{i, \,
    k} \, \left[ \psi^D_{i, \, k}(x) \right]^2
\end{equation}

\end{proposizione}

{\it Proof}: We start with the computation of the sharp of the
function f(x), the definition of the sharp operator (\ref{def:sharpe})
and the chosen structure (\ref{eqn:sharp-indep}) give the following
relation:

\begin{equation}\label{eqn:sharp-terminal-function}
  \begin{array}{lll}
 \displaystyle f(x)^{\#} & = & \displaystyle \sum_k \alpha_k^{\#}
 \, \phi^D_{0, \, k}(x) + \sum_i \sum_k \beta_{i, \,
    k}^{\#} \, \psi^D_{i, \, k}(x) \\
    \scriptstyle & \scriptstyle & \scriptstyle \\
  & = &  \displaystyle \sum_k \sqrt{\gamma(k)} \, \alpha_k \,
   \widehat{\alpha}_k \,  \phi^D_{0, \,
    k}(x) + \sum_i \sum_k \sqrt{\gamma(i, \, k)} \, \beta_{i, \,
    k} \, \widehat{\beta}_{i, \, \, k} \, \psi^D_{i, \, k}(x) \\
    \scriptstyle & \scriptstyle & \scriptstyle
  \end{array}
\end{equation}

Now the second property of operator sharp (see definition
\ref{def:sharpe}) gives the result
(\ref{eqn:gamma-terminal-function}).

The previous theorem as an equivalent for the solution of LPDE.
\begin{flushright}
$\Box$
\end{flushright}

\begin{teorema}[Variance of LPDE solution]

The solution of LPDE has the following variance:

\begin{equation}\label{eqn:variance-solution}
  \Gamma\left[ \widetilde{Q}(t, \, x) \right]  =
  \sum_k \gamma(k) \, \alpha_k^2  \,  \left[\Phi^D_{0, \,
    k}(t, \, x)\right]^2 + \sum_i \sum_k \gamma(i, \, k) \, \beta^2_{i, \,
    k} \, \left[ \Psi^D_{i, \, k}(t, \,x) \right]^2
\end{equation}

\end{teorema}

The proof is similar to the previous theorem, we start with the
sharp decomposition:

\begin{equation}\label{eqn:sharp-solution}
  \begin{array}{lll}
  \displaystyle \widetilde{Q}^{\#}(t, \, x) & = & \displaystyle
   \sum_k \alpha_k^{\#} \, \Phi^D_{0, \, k}(t, \, x) + \sum_i \sum_k \beta_{i, \,
    k}^{\#} \, \Psi^D_{i, \, k}(t, \, x) \\
    \scriptstyle & \scriptstyle & \scriptstyle \\
  & = & \displaystyle  \sum_k \sqrt{\gamma(k)} \, \alpha_k \,
  \widehat{\alpha}_k \,  \Phi^D_{0, \,
    k}(t, \, x) + \sum_i \sum_k \sqrt{\gamma(i, \, k)} \, \beta_{i, \,
    k} \, \widehat{\beta}_{i, \, \, k} \, \Psi^D_{i, \, k}(t, \,x)
    \\
    \scriptstyle & \scriptstyle & \scriptstyle
  \end{array}
\end{equation}

Now the second property of operator sharp, see definition
(\ref{def:sharpe}), gives the result (\ref{eqn:variance-solution}).
\begin{flushright}
$\Box$
\end{flushright}

\begin{remarque}
This uncertainty is easy to estimate due to the fact that of
diffusionlets are universal, in the sense that they are independent
of the initial state, see remark \ref{rem:advantage}.
\end{remarque}

\begin{teorema}[Covariance of LPDE solution]

The solution of LPDE has the following covariance:

\begin{equation}\label{eqn:covariance-solution}
  \begin{array}{lll}
  \displaystyle \Gamma\left[ \widetilde{Q}(t, \, x), \, \widetilde{Q}(s, \, y)  \right] & = &
 \displaystyle \sum_k \gamma(k) \, \alpha_k^2  \,  \Phi^D_{0, \,
    k}(t, \, x) \; \Phi^D_{0, \,
    k}(s, \, y) \\
\scriptstyle & \scriptstyle & \scriptstyle \\
    & & \displaystyle + \sum_i \sum_k \gamma(i, \, k) \, \beta^2_{i, \,
    k} \, \Psi^D_{i, \, k}(t, \,x) \; \Psi^D_{i, \, k}(s, \,y) \\
    \scriptstyle & \scriptstyle & \scriptstyle
     \end{array}
\end{equation}

\end{teorema}

The proof is equal to the previous theorem.
\begin{flushright}
$\Box$
\end{flushright}

\section{Applications in Finance}

In the domain of mathematical finance, the Partial Differential
Equations have a key role. Excepted few cases, the equations in a
financial model do not have a closed form solution, therefore a
numerical approach is mandatory. The Feynman-Kac formula gives us a
bridge between the PDE and SDE, hence it is possible to choose the
numerical method, PDE or Monte-Carlo, depending on the
time-efficiency; generally PDE approach is used in the case of a low
number of variables.

A first application is the Cox-Ingersoll-Ross stochastic
differential equation, see \cite{bib:CIR}:

\begin{displaymath}
dX_t = b\, (a-X_t) \, dt + \sigma \, \sqrt{X_t} dW_t
\end{displaymath}

 when the mean $a$ is zero; the associate PDE is

 \begin{equation*}
 \frac{\partial U}{\partial t} \,(t, \,x) -b \, x \, \frac{\partial U}{\partial
 x} +  \sigma^2 \, x \, \frac{\partial^2 U}{\partial x^2}
\end{equation*}

This case has a relative interest, since an exact characterization
of the solution exists, see Lamberton and Lapeyre
\cite{bib:Lamberton} or Shreve \cite{bib:Shreve}. The principal
advantage consists on the possibility to test the wavelets
procedure.

A second, more relevant, application is a simplified  SABR model;
Hagan et al, see \cite{bib:Hagan}, emphasize the incoherence of the
dynamic behavior of log-normal model, the well-known Black Scholes
model, compared to the behavior observed in the marketplace. In
order to eliminate this problem, Hagan an Woodward, see
\cite{bib:Hagan-woodward}, propose a local volatility model in which
the forward value satisfies

\begin{eqnarray*}
dF_t & =& \sigma_t \, F_t^{\lambda} \, dW_t
\end{eqnarray*}

where $\lambda$ is a fixed parameter, that takes values between 0 and
1, estimated on the market. This model is the starting point for the
SABR model, the SDE is of type (\ref{eqn:LPDE-EDS}); therefore the
procedure described in this article can optimize the procedure of
option pricing.

Hagan et al, see \cite{bib:Hagan}, emphasize an other relevant
aspect, actually market smiles are managed using Dupire's local
volatility models, but a local volatility function different to a
power\footnote{In the case of a power the local volatility model is
reduced at a simplified SABR model.} introduce an intrinsic "length
scale" for the forward price, this inhomogeneous has an hard
financial explanation. Therefore the model proposed by Hagan and
Woodward has the principal problems of all local volatility model,
i.e. a poor equivalent volatility for basket options and an unsmiled
forward volatility, we propose to solve these inconveniences using
the same strategy proposed by Scotti, see \cite{bib:Scotti}, that is
consider a perturbation on the volatility parameter $\sigma_t$ and
handle this uncertainty through an error structure; this study will
be investigate in a following paper.

\section{Conclusion}

In this paper we have investigate the relations between three
objects, i.e. the linear partial differential equations, the error
theory using Dirichlet forms and the wavelets.

These three objects have a different origin and an application
fields separate to date; this article shows as capitalize the
advantages of wavelets bases in order to solve a LPDE and as exploit
the wavelets as a decomposition basis to study the sensibility of
the LPDE.

In a particular case, when the LPDE has the following form:

\begin{displaymath}
U_t = \sigma^2 \, x^{2\, \lambda} U_{xx}
\end{displaymath}

We have proof that the properties of wavelets are partly preserved,
especially the invariance under a scaling and a translation; these
properties permit a fast processing of the general solution of this
type of PDE. The consider LPDE is the new key element of the
non-lognormal financial models, studded by Hagan et al, see
\cite{bib:Hagan-woodward} and \cite{bib:Hagan}.

The principal new feature is the study of the sensibility of the PDE
solution using error theory by means Dirichlet forms, introduced by
Bouleau, see \cite{bib:Bouleau-erreur}, this methodology has our
permitted an evaluation of the sensibility of the solution w.r.t. an
uncertainty on the terminal value, this inaccuracy on the payoff can
model an imprecision on the final spot value, due to the spread
bid-ask as an example.

This paper is just a starting point for the study of inaccuracy on
terminal value and non-lognormal models, in following papers will be
evaluated the impacts of an uncertainty on the volatility parameter
$\sigma$ and on the implied curve of volatility.


\addcontentsline{toc}{section}{Bibliographie}

\begin{thebibliography}{9}

\bibitem{bib:Albeverio} Albeverio, S. (2003): {\it Theory of
    Dirichlet forms and application}, Springer-Verlag, Berlin.

\bibitem{bib:Bouleau-Hirsch} Bouleau, N.; Hirsch, F. (1991):
  {\it Dirichlet Forms and Analysis on Wiener space}, De Gruyter, Berlin.

 \bibitem{bib:Bouleau-erreur} Bouleau, N. (2003): {\it Error
     Calculus for Finance and Physics}, De Gruyter, Berlin.

\bibitem{bib:CIR} Cox, J.C.; Ingersoll, J.E.; Ross, S. (1985): {\it A theory
 of the term structure of interest rates}, Econometrica, 53, pp
 373-384.

   \bibitem{bib:Daubechies-art} Daubechies, I. (1988): {\it Orthonormal Bases of
       Compactly supported wavelets}, Comm. Pure and Appl. Math., 41 pp 909-996.

     \bibitem{bib:Daubechies} Daubechies, I. (1992): {\it Ten Lectures on Waveltes},
      SIAM, Philadelphia.

\bibitem{bib:Daubech-Laga} Daubechies, I.; Lagarias, J. (1992): {\it Two-scale
    difference equations: II. Local Regularity, infinite products of matrices and
    fractals}, SIAM J. of Math. Anal., 24.










\bibitem{bib:Hagan-woodward} Hagan, P.;  Woodward, D.
(1999) {\it Equivalent Black Volatilities} Applied Mathematical
Finance 6, 3, pp. 147-157.


\bibitem{bib:Hagan} Hagan, P.; Kumar, D.; Lesniewsky, A.; Woodward, D.
(2002) {\it Managing Smile Risk} Willmot magazine 1, pp 84-108.

\bibitem{bib:Hardle} Hardle, W.; Kerkyacharian, G.; Picard, D.; Tsybakov, A.
(1998) {\it Wavelets, Approximation and Statistical Applications}, Springer-Verlag, Berlin.

\bibitem{bib:Lamberton} Lamberton, D.; Lapeyre, B. (1995):
{\it Introduction to Stochastic Calculus Applied to Finance}, Chapman \& Hall, London.


\bibitem{bib:Mallat} Mallat, S.
(2000) {\it Une exploration des signaux en ondelettes}, Editions de l'Ecole Polytechnique, Paris.


\bibitem{bib:FWT} Mallat, S. (1989) {\it A theory for multiresolution signal
decomposition: the wavelet representation }, IEEE Pattern Anal. and
Machine Intell., vol. 11, no. 7, pp. 674-693

\bibitem{bib:Scotti} Scotti, S.;  {\it Perturbative Approach on Financial Markets}, submitted to Mathematical Finance.


\bibitem{bib:Shen-Strang} Shen, J.; Strang, G.
(2000) {\it On Wavelet Fundamental Solutions to the Heat Equation -
Heatlets}, Journal of Differential Equations, vol. 161, no. 2,
pp.403-421 .

\bibitem{bib:Shreve} Shreve, S.E. (2004) {\it Stochastic Calculus
for Finance II, Continuous-Time Models}, Springer, New York







\end{thebibliography}
\end{document}